\documentclass[12pt]{article}

\newif\ifpdf
\ifx\pdfoutput\undefined
    \pdffalse       
\else
    \pdfoutput=1    
    \pdftrue
\fi

\usepackage{amsmath,amssymb,amsthm}

\ifpdf
    \usepackage[pdftex]{graphicx}
    \usepackage[backref, colorlinks=true, pdfstartview=FitV, linkcolor=blue,
        citecolor=blue, urlcolor=blue]{hyperref}
\else
    \usepackage{graphicx}
    \usepackage{hyperref}
\fi
\usepackage[numeric]{amsrefs} 

    \font\smallit=cmti10
    \font\smalltt=cmtt10
    \font\smallrm=cmr9

\title{On Z.-W. Sun's Disjoint Congruence Classes Conjecture}
\author{Kevin O'Bryant \\ City University of New York (College of Staten Island)\\ \url{kevin@member.ams.org}}
\date{\today}

    \newcommand{\ceiling}[1]{\mbox{$\big\lceil #1 \big\rceil$}}

    \DeclareMathOperator{\lcm}{LCM}

    \newtheorem{thm}{Theorem}
    \newtheorem{lem}[thm]{Lemma}

    \newtheorem{cnj}[thm]{Conjecture}
    \newtheorem{prop}[thm]{Proposition}

\begin{document}
    \ifpdf
       \DeclareGraphicsExtensions{.pdf,.jpg,.mps,.png}
    \fi
\begin{center}
 {\bf ON Z.-W. SUN'S DISJOINT CONGRUENCE CLASSES CONJECTURE}
 \vskip 20pt
 {\bf Kevin O'Bryant}\\
 {\smallit Department of Mathematics, City University of New York, \\College of Staten Island, Staten Island, NY
 10314, U. S. A.}\\
 {\url{kevin@member.ams.org}}\\
\end{center}
 \vskip 30pt
\centerline{\smallit Received: , Accepted: , Published: }
 \vskip 30pt

\centerline{\bf Abstract}

\noindent Sun has conjectured that if $k$ congruence classes are disjoint, then necessarily two of the moduli
have greatest common divisor at least as large as $k$. We prove this conjecture for $k\le 20$.

\pagestyle{myheadings} \markright{\smalltt INTEGERS: \smallrm ELECTRONIC JOURNAL OF COMBINATORIAL NUMBER THEORY
\smalltt x (200x), \#Axx\hfill}

\thispagestyle{empty} \baselineskip=15pt \vskip 30pt

\section*{\normalsize 1. Introduction} In May of 2003, Prof Sun\footnote{Zhi-Wei Sun, Department of Mathematics, Nanjing University,
\url{zwsun@nju.edu.cn}} proposed a conjecture on the number theory
listserver\footnote{\url{http://listserv.nodak.edu/archives/nmbrthry.html}}.

\begin{cnj}[Disjoint Congruence Classes Conjecture]\label{cnj:dccc}
    If $k\ge 2$ congruence classes $a_i \pmod{m_i}$ are disjoint, then there exist $i<j$ with
    $\gcd(m_i,m_j)\ge k$.
\end{cnj}

The set of congruence classes $\{1\pmod{k},2\pmod k, \dots, k \pmod k\}$ demonstrate that, if true, the DCCC is
best possible.

The contrapositive of the $k=2$ case is a familiar special case of the Chinese remainder theorem: if two
congruence classes have relatively prime moduli, then they intersect. Prof Graham\footnote{Ronald Graham,
Department of Computer Science and Engineering, University of California, San Diego, \url{graham@ucsd.edu}}, who
first brought this problem to my attention, pointed out that the $k=3$ case follows easily from the pigeonhole
principle, as we now explain.

Suppose that $a_1 \bmod{m_1}$, $a_2\bmod{m_2}$, $a_3\bmod{m_3}$ is a counter-example: they are disjoint and
$\gcd(m_i,m_j)<3$ for all $i<j$. If $\gcd(m_i,m_j)=1$, then by the Chinese remainder theorem $a_i\pmod{m_i}$ and
$a_j\pmod{m_j}$ intersect. Thus, $\gcd(m_i,m_j)=2$ for all $i<j$; in particular, all $m_i$ are even. By the
pigeonhole principle, two of the $a_i$'s must have the same parity, say $a_1\equiv a_2\pmod{2}$. Since
$\gcd(m_1,m_2)=2$, obviously $\gcd(\frac{m_1}{2},\frac{m_2}{2})=1$. If $a_1$ and $a_2$ are both even, then by
the Chinese remainder theorem we can find $x$ in the intersection of
    \[
    \frac{a_1}{2} \pmod{\frac{m_1}{2}} \quad\text{ and }\quad  \frac{a_2}{2} \pmod{\frac{m_2}{2}}
    \]
and consequently $2x$ is in the intersection of $a_1 \pmod {m_1}$ and $a_2 \pmod{m_2}$. If $a_1$ and $a_2$ are
both odd, then we can find $x$ in the intersection of
    \[
    \frac{a_1-1}{2} \pmod{\frac{m_1}{2}} \quad \text{ and }\quad \frac{a_2-1}{2} \pmod{\frac{m_2}{2}}.
    \]
But then $2x+1$ is in the intersection of $a_1 \pmod {m_1}$ and $a_2 \pmod{m_2}$.

The following criterion is at the heart of the last paragraph.
\begin{prop}\label{Fraenkel}
The congruence classes $a_i\pmod{m_i}$ and $a_j\pmod{m_j}$ are disjoint if and only if $\gcd(m_i,m_j) \nmid
a_i-a_j$.
\end{prop}

Prof Graham also noted that the $k=4$ case was similar but with more considerations, and considered it likely
that $k=5$ was similarly tractable. We now state our main theorem.

\begin{thm}
The DCCC holds for $k\le 20$. Moreover, a counterexample to the DCCC with minimal $k$ does not have
$k\in\{24,30\}$.
\end{thm}

We close the introduction with a stronger conjecture of Prof Sun~\cite{Sun2}.
\begin{cnj}
Suppose that $A_1, \dots, A_k$ are disjoint left-cosets of $H_1,\dots, H_k$ in the group $G$. Then
$\gcd([G:H_i],[G:H_j]) \ge k$ for some $i<j$.
\end{cnj}

\vskip 30pt
\section*{\normalsize 2. A disjointness criterion}
The following criterion is stated in~\cite{Huhn.Megyesi} without proof.
\begin{lem}\label{lem:test}
If
    $\displaystyle \sum_{i=1}^\ell \frac{1}{\gcd(m_i,M)} > 1,$
where $M$ is any multiple of $\,\lcm\{\gcd(m_i,m_j) \colon 1\le i < j\le \ell\}$, then the congruence classes
$a_i\bmod{m_i}$ (with $1\le i \le k$) are not disjoint.
\end{lem}

The usefulness of this criterion lies in the fact that it makes no reference to the $a_i$.

Can there be seven disjoint congruence classes with moduli 20, 15, 12, 6, 6, 6, 6? Lemma~\ref{lem:test} does not
exclude it directly, but it does exclude the possibility that there are six disjoint congruence classes with
moduli 15, 12, 6, 6, 6, 6. It is plausible that if every subset of $\{m_1,\dots,m_k\}$ passes the above test,
then it is possible to choose $a_1,\dots,a_k$ so that $a_1\bmod{m_1},\dots,a_k\bmod{m_k}$ are disjoint, and
Huhn~\&~Megyesi~\cite{Huhn.Megyesi} conjectured this. However, Z.-W. Sun~\cite{Sun1} notes that the moduli 10,
15, 36, 42, 66 pass the above test (as does every subset of the moduli), but these are not the moduli of
disjoint congruence classes.

\begin{proof}
The class $a_i\bmod m_i$ intersects each of the $M/\gcd(m_i,M)$ classes $a_i + j m_i \pmod M$ (for $0\le j <
{M}/{\gcd(m_i,M)}$); that is, the class $a_i \bmod m_i$ is actually $M/\gcd(m_i,M)$ classes modulo $M$. Since by
hypothesis $\sum_{i=1}^\ell M/\gcd(m_i,M) >M$, the pigeonhole principle implies that two of the modulo-$M$
congruence classes must intersect. In other words, there are integers $\alpha,\beta,\gamma,i,j$ such that
    \[a_i+\alpha m_i = a_j + \beta m_j +\gamma M.\]

Since $\gcd(m_i,m_j) \mid M$, there are integers $\delta,\epsilon$ such that $\delta m_i+\epsilon m_j=M$. Thus
    \[a_i+(\alpha-\delta\gamma) m_i = a_j +(\beta+\epsilon\gamma) m_j,\]
and we see that $a_i \bmod m_i$ and $a_j \bmod m_j$ intersect.
\end{proof}

\vskip 30pt
\section*{\normalsize 3. Without Loss of Generality}
We suppose that a minimal counterexample exists, and use that to determine a sequence $m_1,\dots,m_k$ (as
described in the lemma below) that has particular properties (in particular, we have an explicit upper bound on
$m_i$). While there are infinitely many sets of $k$ congruence classes, there are only finitely many such $m_i$
sequences, and in fact we show that for $k\leq 20$ there are none. We note that the existence of such an $m_i$
sequence would not disprove Sun's Disjoint Congruence Classes Conjecture, but that nonexistence would imply his
conjecture.

Set $L_{k}:= \lcm\{1,2,\dots,k-1\}$.
\[ \begin{array}{c|ccccccccccccc}
    k   & 3 & 4 & 5 & 6 & 7 & 8 & 9 & 10 & 11 & 12 & 13 & 14 & 15\\\hline
    L_k & 2 & 6 &12 &60 &60 &420& 840& 2520&2520& 27720 & 27720 & 360360 & 360360 \\
    \end{array}
\]

\begin{lem}\label{lem:G}
Suppose that $k\ge 2$ is the least integer such that there are $k$ disjoint congruence classes $a_i\pmod{m_i}$
with $\gcd(m_i,m_j)<k$ for $i<j$ (i.e., a counterexample to the DCCC). Further, suppose that $a_i\pmod{m_i}$
($1\le i \le k$) has minimal $\sum_i m_i$. Set
    \[N:=\lcm\{\gcd(m_i,m_j) \colon 1\le i<j\le k\}. \]
Then $k\ge 4$, and for all $i$
    \begin{enumerate}
        \item   $m_i \mid N$, and $N \mid L_k$, and $N=\lcm\{m_1,m_2,\dots,m_k\}$;\label{item:Lk}
        \item   $1<\gcd(m_i,m_j)<k$ for all $j\not=i$;\label{item:gcd}
        \item   If a prime power $q$ divides $m_i$, then there is $j\not=i$ with $q \mid m_j$.\label{item:share}
        \item   $m_i$ is not a prime power;\label{item:primepower}
        \item   At least three of the $m$'s are multiples of $k-1$;\label{item:minimal}
        \item   If only three of the $m$'s are multiples of $k-1$, then two others are multiples of
        $k-2$;\label{item:minimal2}
        \item   For every $\{h_1,\dots,h_\ell\}\subseteq\{m_1,\dots,m_k\}$, and every multiple $M$ of \newline
        \mbox{$\lcm\{\gcd(h_i,h_j)\colon 1\le i<j\le\ell\}$},
            \[\sum_{i=1}^\ell \frac1{\gcd(h_i,M)} \le 1.\]  \label{item:test}
        \item   If $7\le k \le 30$, and $p\ge k/2$ is a prime that divides some $m$, then it divides exactly two
        of the $m$'s.\label{item:bigprime}
    \end{enumerate}
\end{lem}

\begin{proof}
We noted in the introduction the reasons why $k\not\leq 3$.

In this paragraph, we prove the three statements in item~\ref{item:Lk}. Suppose that $p^r$ is a prime power that
divides $m_i$ but not $N$. By the definition of $N$, we see that $p^r \nmid m_j$ for $j\not=i$. Thus
$\gcd(m_i,m_j)=\gcd(m_i/p,m_j)$, and by Proposition~\ref{Fraenkel}, the classes $a_i \pmod{m_i/p}$ and
$a_j\pmod{m_j}$ are disjoint. Thus, we can replace $a_i \pmod{m_i}$ in our counterexample to the DCCC with $a_i
\pmod{m_i/p}$, obtaining a counterexample with smaller $\sum m$. We began with minimal $\sum m$, so we conclude
that there is no prime power dividing $m_i$ but not $N$. Since $N$ is defined to be the least common multiple of
a subset of $\{1,2,\dots,k-1\}$, it is clear that $N\mid L_k$. Moreover, since $N$ is the least common multiple
of divisors of the $m_i$, and the $m_i$ divide $N$, it is also now immediate that $N=\lcm\{m_1,\dots,m_k\}$.

That $\gcd(m_i,m_j)<k$ is by hypothesis; that $\gcd(m_i,m_j)>1$ follows from the chinese remainder theorem. This
proves item~\ref{item:gcd}.

Any prime power that divides some $m_i$ also divides $N$ since $N=\lcm\{m_1,\dots,m_k\}$. But any prime power
that divides $N=\lcm\{\gcd(m_i,m_j)\colon i\not= j\}$ must divide two of the $m$'s. This proves
item~\ref{item:share}.

Items~\ref{item:primepower},~\ref{item:minimal}, and~\ref{item:minimal2} are based on the minimality of $k$.
Suppose that $m_1=p^r$, a prime power. Since $\gcd(m_1,m_j)>1$, we know that $p\mid m_j$ for every $j$. By
item~\ref{item:share}, there is $m_2$ that is also a multiple of $p^r$. Since $\gcd(m_1,m_2)<k$, we see that
$p<k$, whence $\ceiling{k/p}\ge 2$. By Proposition~\ref{Fraenkel}, $\gcd(m_i,m_j) \nmid a_j-a_i$ for every $1\le
i < j\le k$. And by the pigeonhole principle, there are at least $\ceiling{k/p}$ of the $a$'s that are congruent
to one another modulo $p$, say (after renumbering)
    \[a_1\equiv a_2 \equiv \dots \equiv a_{\lceil k/p \rceil} \pmod{p}.\]
Consequently,
    \[\gcd\left(\dfrac{m_i}p,\dfrac{m_j}p\right) \nmid \frac{a_j-a_i}{p}\]
for every $1\le i < j\le {\lceil k/p \rceil}$. This implies, again by Proposition~\ref{Fraenkel}, that
    \[0 \pmod{\dfrac{m_1}p},\quad \frac{a_2-a_1}{p} \pmod{\dfrac{m_2}p},\quad \dots,\quad
        \frac{a_{\lceil k/p \rceil}-a_1}{p} \pmod{\frac{m_{\lceil k/p \rceil}}p}\]
are also disjoint. Since $k$ is minimal, there must be $1\le i < j \le {\lceil k/p \rceil}$ with
    \[\gcd \left( \frac{m_i}{p}, \frac{m_j}{p} \right) \geq {\left\lceil \frac{k}{p} \right\rceil},\]
whence $\gcd(m_i,m_j) \ge {\lceil k/p \rceil} p \ge k$, contradicting the existence of a counterexample with $k$
classes. This proves item~\ref{item:primepower}.

Suppose that only $m_1$ and $m_2$ are multiples of $k-1$. Then $a_i \pmod{m_i}$ (with $2\le i \le k$) is a
collection of $k-1$ disjoint congruence classes whose moduli have gcds strictly less than $k-1$. This
contradicts the minimality of $k$, and proves item~\ref{item:minimal}.

Suppose that only $m_1, m_2,$ and $m_3$ are multiples of $k-1$. Then $a_i\pmod{m_i}$ (with $3\le i \le k$) are
$k-2$ disjoint congruence classes whose moduli have gcds strictly less than $k$ (because we started with a
counterexample). Moreover, the gcds are strictly less than $k-1$ because none of $m_4,\dots,m_k$ are multiples
of $k-1$. If there are not two of $m_3, \dots, m_k$ that are multiples of $k-2$, then we have even more: a
counterexample with $k-2$ sequences. By the minimality of $k$, then, two of $m_3,\dots,m_k$ must be multiples of
$k-2$. If there are exactly two, then it could be that $m_3$ is one of them. To see that is not the case, we
renumber $m_2$ and $m_3$ to conclude that both $m_2$ and $m_3$ are multiples of $k-2$. But then both $k-1$ and
$k-2$ divide $\gcd(m_2,m_3)$, so that it must be at least $k$. This proves item~\ref{item:minimal2}, and that
$k\ge 5$.

Item~\ref{item:test} is a restatement of Lemma~\ref{lem:test} for subsets.

Now suppose that $p$ and $k$ are as hypothesized in item~\ref{item:bigprime}, and suppose that
$m_1,\dots,m_\ell$ are multiples of $p$, while $m_{\ell+1},\dots,m_k$ are not multiples of $p$. Note that
$\ell\ge2$ by item~\ref{item:share}, and by item \ref{item:gcd} necessarily $p<k$. Suppose that there are $r$
primes less than $k$; since $k\le 30$ we know that $r\le 10$. Assume, by way of contradiction, that $\ell\ge 3$.

Set
    \[
    {\cal P}_i := \{q \colon q\text{ prime}, q \mid m_i\}\setminus\{p\},\]
and recall that by Lemma~\ref{lem:G} (item \ref{item:primepower}) the cardinality of ${\cal P}_i$ is at least 1
for every $i\le \ell$, and at least 2 for $i>\ell$. Moreover, since $\gcd(m_i,m_j)<k\le 2p$ the ${\cal P}_i$
(for $1\le i \le \ell$) are disjoint.

Now define for $\ell<i\le k$ the $\ell$-tuple
    \[
    \omega_i := \bigg( \min\left\{ {\cal P}_1 \cap {\cal P}_i\},\min\{{\cal P}_2 \cap {\cal P}_i\}, \dots,
                        \min\{{\cal P}_\ell \cap {\cal P}_i\}\right\} \bigg),
    \]
which is an element of ${\cal P}_1 \times \cdots \times {\cal P}_\ell$. If $\omega_i=\omega_j$ (with $\ell < i
<j\le k$), then $\gcd(m_i,m_j)$ is at least the product of the primes in $\omega_i$; since $\ell\ge 3$ this
product must be at least $2\cdot 3\cdot 5=30\ge k$. By the pigeonhole principle, this certainly happens if
    \[
    k-\ell> \big|{\cal P}_1 \times \cdots \times {\cal P}_\ell\big| = \prod_{i=1}^\ell |{\cal P}_i|.
    \]
Since $|{\cal P}_i|\ge 1$ and $\sum_{i=1}^\ell |{\cal P}_i| \le r-1$ (because the ${\cal P}_i$ are disjoint for
$i\le \ell$, and there are $r$ primes less than $k$ including $p$), we can easily bound the size of
$\prod_{i=1}^\ell |{\cal P}_i|$ in terms of $r$ and $\ell$:
    \begin{center}
    \hspace{6.6cm} $\ell$ \newline
    \[ r\quad
    \begin{array}{c|cccccccc}
    \max\{\prod_{i=1}^\ell |{\cal P}_i|\} & 3 & 4 & 5 & 6 & 7 & 8 & 9 \\ \hline
    4 & 1 & 0& 0& 0& 0& 0& 0\\
    5 & 2 & 1 & 0& 0& 0& 0& 0 \\
    6 & 4 & 2 & 1 & 0& 0& 0& 0\\
    7 & 8 & 4 & 2 & 1 & 0& 0& 0 \\
    8 & 12 & 8 & 4 & 2 & 1 & 0 & 0 \\
    9 & 18 & 16 & 8 & 4 & 2 & 1 & 0 \\
    10 & 27 & 24 & 16 & 8 & 4 & 2 & 1
    \end{array}
    \]
    \end{center}
This table shows that for $4\le r\le 10$, necessarily $k-\ell>|\text{Range}(\omega)|$, and so two $m$'s have gcd
at least $30\ge k$, with the exception of $r=10$, $k=30$, $\ell=3$.

In the case of $r=10$, $k=30$, $\ell=3$, each of the 27 possible values of $\omega$ must actually occur: $|{\cal
P}_1|=|{\cal P}_2|=|{\cal P}_3| =3$. Thus two of the four primes $17$, $19$, $23$, $29$ are in separate ${\cal
P}_i$'s ($1\le i \le 3$), and so must both be in three of the $\omega_i$'s ($4\le i \le 30$). But then the two
corresponding $m$'s will have gcd at least $17\cdot 19>30$.
\end{proof}

\vskip 30pt
\section*{\normalsize 4. Casework}

\subsection*{\normalsize 4.1. The cases $3\le k \le 6$}
These cases can be handled in a variety of ways. We handle them here using only \mbox{$1<m_i \mid L_k$},
$1<\gcd(m_i,m_j)<k$, and that $m_i$ cannot be a prime power.

\paragraph{{$\mathbf{k=3}.$}}
The divisors of $L_3=2$ are 1 and 2. That $m_i=1$ is impossible since $\gcd(m_i,m_j)>1$, and that $m_i=2$ is
impossible since $m_i$ cannot be a prime power.

\paragraph{$\mathbf{k=4}.$}
The divisors of $L_4=6$ are 1, 2, 3, and 6. Since $m_i$ cannot be a prime power (or 1), every $m_i=6$. But
$\gcd(m_i,m_j)<k<6$, so this too is impossible.

\paragraph{$\mathbf{k=5}.$}
The divisors of $L_5=12$ are 1, 2, 3, 4, 6, and 12. Since $m_i$ cannot be a prime power (or 1), every $m_i$ is
either $6$ or 12. But $\gcd(m_i,m_j)<k<6$, so this too is impossible.

\paragraph{$\mathbf{k=6}.$}
The only allowed divisors of $L_6=60$ are $6, 10, 12, 15, 20, 30,$ and 60. Since $\gcd(m_i,m_j)<6$, the $m_i$
must be distinct. By the pigeonhole principle, two of the $m$'s must be multiples of 10, whence
$\gcd(m_i,m_j)\ge 10$.

\subsection*{\normalsize 4.2. The cases $\mathbf{6\le k \le 10}$}

While looking through the cases below, I find it helpful to imagine the complete graph on $k$ vertices with
vertices labeled with $m_1,m_2,\dots,m_k$, and edges labeled with $\gcd(m_i,m_j)$.

With $k=7$, the assignments $m_1=20$, $m_2=15$, $m_3=12$, $m_4=m_5=m_6=m_7=6$ satisfy all of the conditions
given in Lemma~\ref{lem:G} except item~\ref{item:test}, and also passes the test of Lemma~\ref{lem:test}. Thus,
for $k\ge7$ the arguments are necessarily more involved.

\paragraph{$\mathbf{k=7}.$}
We have $N \mid L_7 = 60$. Suppose that none of the $m$'s are multiples of 5, so that we can take $M:=12$, and
at most one $\gcd(m_i,M)$ is 12, with other six $\gcd(m_i,M)$ being $\le 6$. In this case,
    \[\sum_{i=1}^7 \frac{1}{\gcd(m_i,12)} \ge \frac 1{12}+6\,\frac16 >1.\]

Now suppose that some $m$ is a multiple of 5, and by item~\ref{item:share}, at least two of the $m$'s are
multiples of 5. We know that no $m$ is exactly 5 (a prime), so that the two (or more) of the $m$'s that are
multiples of 5 are in $\{10, 15, 20, 30, 60\}$. Since $\gcd(m_i,m_j)\le 6$ the only possibilities are ``10 and
15'' or ``15 and 20''. Delete the congruence class that gives $m$ being either 10 or 20, and the six classes
remaining fail to have the condition given in item~\ref{item:test} with $M=12$:
    \[\sum_{\text{6 classes}} \frac{1}{\gcd(m_i,12)} \ge \frac{1}{\gcd(15,12)}+\frac{1}{\gcd(12,12)}+
        4\,\frac{1}{\gcd(6,12)} >1.\]

\paragraph{$\mathbf{k=8}.$}
Here, items \ref{item:minimal} and \ref{item:bigprime} are not consistent. For clarity, we reprove here this
special case of item~\ref{item:bigprime}, proved in full above.

By item~\ref{item:minimal}, we may assume that $m_1,m_2,m_3$ are multiples of 7 and divisors of
$L_8=2^2\cdot3\cdot5\cdot7$. No $m$ is exactly 7, so the remaining possibilities for $(m_1,m_2,m_3)$ are
$(2\cdot7, 3\cdot 7, 5 \cdot 7)$ and $(4\cdot7, 3\cdot 7, 5 \cdot 7)$. In particular, there are not {\em four}
$m$'s that are multiples of 7. Item~\ref{item:minimal2} now implies that there {\em are} two $m$'s, say $m_4$
and $m_5$, that are multiples of $k-2=6$, and one of those is not a multiple of 5 (or else $\gcd(m_4,m_5)=30$).
Thus, either $\gcd(m_4,m_3)$ or $\gcd(m_5,m_3)$ is 1.

\paragraph{$\mathbf{k=9}.$}
At least three of the $m$'s are multiples of 8 and divisors of $L_9=2^3\cdot 3\cdot5\cdot7$. Since $m_i\not=8$
(a prime power), the three multiples of 8 are $3\cdot 8, 5\cdot 8, 7\cdot 8$, and there isn't a fourth. Since
there isn't a fourth multiple of 8, there must be two other $m$'s that are multiples of $k-2=7$, one of which is
relatively prime to 5, and both are odd since $\gcd(m_i,7\cdot 8) \le 8$. Thus there is one that is relatively
prime to $5\cdot 8$, a contradiction.

\paragraph{$\mathbf{k=10}.$}
At least three of the $m$'s are multiples of 9 and divisors of $L_{10}=2^3\cdot 3^2\cdot5\cdot7$. Since
$m_i\not=9$ (a prime power) the three multiples of 9 are $(m_1,m_2,m_3)=(2^r\cdot 9, 5\cdot 9, 7\cdot 9)$ (for
some $r\ge1$), and there isn't a fourth $m$ that is a multiple of 9. Since there isn't a fourth multiple of 9,
there must be two other $m$'s, say $m_4$ and $m_5$, that are multiples of 8, one of which is relatively prime to
5, and so a multiple of 3 since $\gcd(m_i,m_5)>1$. Say $m_4$ is the multiple of $3\cdot8$, and $m_5$ is not a
multiple of 3 since $\gcd(m_4,m_5)<10$. Since $\gcd(m_2,m_5)>1$ and $\gcd(m_3,m_5)>1$, we have $m_5=8\cdot 5
\cdot 7$. Since $\gcd(m_1,m_4)<10$ and $\gcd(m_4,m_5)=8$, we have $m_4=3\cdot8$. And since $m_1\not=9$ and
$\gcd(m_1,m_4)<10$, we find that $m_1=2^r \cdot 9 = 2\cdot 9$.

We have $m_1=2\cdot 3^2, m_2=5\cdot 3^2, m_3=7\cdot 3^2, m_4=2^3\cdot 3, m_5=2^3\cdot5\cdot7$. Suppose that
there is another $m$ that is a multiple of 7. Since $\gcd(m,m_1)>1$, either $2 \mid m$ or $3 \mid m$. But then
either $\gcd(m,m_5)=14\ge k$ or $\gcd(m,m_3)=21 \ge k$. Likewise, no other $m$ is a multiple of 5. Thus the
other $m$'s are divisors of $2^3\cdot 3^2$, but not multiples of 8 or 9 or prime powers: the possibilities are
4, 6, 12. No other $m$ can be a multiple of 12 because $m_4=24$, so $m_6,\cdots, m_{10}$ are at most 6. Deleting
$m_5$, we take $M=2^3\cdot3^2$ and
    \[\sum_{i\not=5} \frac{1}{\gcd(m_i,M)} \geq \frac 1{18}+\frac 19 + \frac 19 + \frac 1{24}+\frac 18+5\,\frac
    16 >1.\]

\subsection*{\normalsize 4.3. The cases $\mathbf{k\in\{8,12,14,18,20,24,30\}}$}
For $k\in\{8,12,14,18,20,24,30\}$, the prime $k-1$ must divide at least 3 of the $m$'s by
item~\ref{item:minimal}. By item~\ref{item:bigprime}, however, only 0 or 2 of the $m$'s can be multiples of a
prime as large as $k/2$.

\subsection*{\normalsize 4.4. The cases $\mathbf{k \le 19}$}

Let {\tt G} be a (possibly empty) list of positive integers and let {\tt potentials} be a (possibly empty) list
of positive integers. Define a function {\tt Grow[G,potentials]} that will return {\tt True} if it is possible
to augment {\tt G} with elements of {\tt potentials} to get a set of size $k$ (a global variable) such that each
pair of elements has $\gcd$ strictly between 1 and $k$, and every subset of it passes the test of
Lemma~\ref{lem:test}. Otherwise it will return {\tt False}. We give in Figure~\ref{code} a {\em Mathematica}
program that accomplishes this.

\begin{figure}[b]
    {\tt
    \begin{verbatim}
    GCDList[G_] := Union[Flatten[Table[GCD[G[[i]], G[[j]]],
                                       {i, 1, Length[G]-1},
                                       {j, i+1, Length[G]}    ]]];
    LemmaTest[G_]:= Block[
                          {gcdlist = GCDList[G],
                           M },
                          M = LCM @@ gcdlist;
                          ( (Last[gcdlist] < k) &&
                            (First[gcdlist] > 1) &&
                            (Plus @@ (1/ (GCD[M,#]& /@ G) ) <= 1)
                          )                                   ];
    Grow[G_, potentials_] :=
        If[Length[G] == k,
           {passed, counter} = {True, Length[Subsets[G, {3, k}]]};
           While[passed && counter > 0,
                 passed = LemmaTest @@ Subsets[G, {3, k}, {counter}];
                 counter--                                    ];
           passed,
           Or @@ Table[Grow[Append[ G, potentials[[i]] ],
                            Select[potentials,
                                   ( (# >= potentials[[i]]) &&
                                     LemmaTest[Join[G, {potentials[[i]], #}]]
                                   )&                         ]],
                       {i, Length[potentials]}                ]];
    \end{verbatim}
    }
    \caption{{\em Mathematica} code to look for counter-examples.\label{code}}
\end{figure}

We can prove that there is no counter-example with $k\leq 19$ sequences by executing the following {\em
Mathematica} loop
    {\tt \begin{verbatim}
    Or @@ Table[L = LCM @@ Range[k-1];
                candidates = Select[Divisors[L],
                                    (Length[FactorInteger[#]] > 1)&];
                Grow[ {}, candidates ],
                {k, 2, 19}                                    ]
    \end{verbatim}}
\noindent and observing that the output is {\tt False}. This required slightly less than a week to run on the
authors humble desktop PC.

\vskip 30pt
\section*{\normalsize Acknowledgement} The author wishes to thank Prof Landman\footnote{Bruce M. Landman,
Department of Mathematics, State University of West Georgia, \url{landman@westga.edu}} for reproducibly
marvelously organizing the {\em Integers} conference. He also wishes to thank Prof Eichhorn\footnote{Dennis
Eichhorn, Department of Math and Computer Science, California State University, East Bay,
\url{eichhorn@mcs.csueastbay.edu}} for several conversations concerning the DCCC, including a discussion that
lead to Lemma 6 (item 8).

\end{document}